\documentclass[runningheads,a4paper]{llncs}
\usepackage{amssymb}
\usepackage{graphicx}
\usepackage{url}
\newcommand{\keywords}[1]{\par\addvspace\baselineskip
\noindent\keywordname\enspace\ignorespaces#1}

\usepackage{bbm}
\usepackage{amsmath}
\usepackage{amsfonts}
\usepackage{color}
\usepackage{url}
\usepackage{microtype}
\usepackage{adjustbox}
\usepackage{tikz}
\usepackage{hyperref}
\usepackage{xspace}
\definecolor{darkblue}{RGB}{0,0,160}
\hypersetup{
colorlinks,%
citecolor=black,%
filecolor=black,%
linkcolor=darkblue,%
urlcolor=darkblue
}

\usepackage[utf8]{inputenc}

\newcommand{\excise}[1]{}%{$\star$\textsc{#1}$\star$}

\newcommand{\ring}[1]{\ensuremath{\mathbb{#1}}}

\newcommand\RR{\ring{R}}

 \DeclareMathOperator{\supp}{supp} % The support
 \DeclareMathOperator{\diag}{diag} % The diagnol matrix
 \DeclareMathOperator{\conv}{conv} % convex hull

\newcommand{\yalmip}{\textsc{yalmip}\xspace}
\newcommand{\mathematica}{\textsc{Mathematica}\xspace}
\newcommand{\qepcad}{\textsc{qepcad}\xspace}

\newcommand{\mosek}{\textsc{mosek}\xspace}

\begin{document}

\mainmatter
%================================================

%==== FILL IN ====================================
\title{On the feasibility of semi-algebraic sets in Poisson regression}  % Full title
\titlerunning{Semi-algebraic sets in Poisson regression} % Shor ttitle
\author{Thomas Kahle}
\authorrunning{Kahle}
\institute{
Otto-von-Guericke Universität Magdeburg, Germany\\
\email{thomas.kahle@ovgu.de},\\ 
\texttt{https://www.thomas-kahle.de}
}
\maketitle

\begin{abstract}
Designing experiments for generalized linear models is difficult
because optimal designs depend on unknown parameters.  The local
optimality approach is to study the regions in parameter space where a
given design is optimal.  In many situations these regions are
semi-algebraic.  We investigate regions of optimality using computer
tools such as \yalmip, \qepcad, and \mathematica.
% 
% Here we
% investigate local optimality.  We propose to study for a given design
% its region of optimality in parameter space.  Often these regions are
% semi-algebraic and feature interesting symmetries.  We demonstrate
% this with the Rasch Poisson counts model.  For any given interaction
% order between the explanatory variables we give a characterization of
% the regions of optimality of a special saturated design. This extends
% known results from the case of no interaction.  We also give an
% algebraic and geometric perspective on optimality of experimental
% designs for the Rasch Poisson counts model using polyhedral and
% spectrahedral geometry.
\keywords{algebraic statistics, optimal experimental design, Poisson
regression, semi-algebraic sets}
\end{abstract}

\section{Introduction}
\label{sec:introduction-1}

Generalized linear models are a mainstay of statistics, but optimal
experimental designs for them are hard to find, as they depend on
unknown parameters of the model.  A common approach to this problem is
to study local optimality, that is, determine an optimal design for
each fixed set of parameters.  In practice, this means that
appropriate parameters have to be guessed a priori, or fixed by other
means.  In \cite{kahle2015algebraic} the authors approached this
problem from a global perspective.  They study the \emph{regions of
optimality} of fixed designs and demonstrate that these are often
defined by semi-algebraic constraints.  Their main tool is a general
equivalence theorem due to Kiefer and Wolfowitz, which directly yields
polynomial inequalities in the parameters.  This makes these problems
amenable to the toolbox of real algebraic geometry.  In this extended
abstract we pursue this direction for the Rasch Poisson counts model
which is used in psychometry~\cite{doebler2015processing} in the
design of mental speed tests.  Analyzing saturated designs for this
model amounts to studying the feasibility of polynomial inequality
systems.  We examine the state of computer algebra tools for this
purpose and find that there is room for improvement.

\subsection*{Acknowledgement}
The author is supported by the Research Focus Dynamical Systems (CDS)
of the state Saxony-Anhalt.

\section{Polynomial inequality systems in statistics}

For brevity we omit any details of statistical theory and focus on
mathematical and computational problems.  The interested reader should
consult \cite{kahle2015algebraic} and its references.  We also stick
to that paper's notation.
Throughout, fix a positive integer $k$, the \emph{number of rules},
and another positive integer $d\le k$, the \emph{interaction order}.
A \emph{rule setting} is a binary string
$x = (x_1,\dots,x_k) \in\{0,1\}^k$.
% We freely identify $x$ with the
% corresponding subset of $[k]:=\{1,\dots,k\}$ and write, for example,
% $x\subseteq A$ if $A \subset [k]$.  
The \emph{regression function of interaction order $d$} is the
function $f:\{0,1\}^k \to \{0,1\}^p$ whose components are all
square-free monomials of degree at most~$d$ in the indeterminates
$x_1,\dots,x_k$.
The value $p$ equals the number of square-free monomials of degree at
most $d$ and depends on $d$ and~$k$.  
% The regression function takes values in $\{0,1\}^p$, but we use
% $\RR^p$ as its codomain so that the following is well defined.
For any $\beta\in\RR^p$, the \emph{intensity} of the rule setting
$x\in\{0,1\}^k$ is
\[
\lambda(x,\beta) = e^{f(x)^T\beta}.
\]
The \emph{information matrix} of $x$ at $\beta$ is the rank one matrix
\[
M(x,\beta) = \lambda(\beta,x) f(x)f(x)^T.
\]
The \emph{information matrix polytope} is
\[
P(\beta) = \conv \{M(x,\beta) : x\in\{0,1\}^k\}.
\]
The case $d=1$ and $k$ arbitrary is known as the model with \emph{$k$
independent rules}.  In this case $f(x) = (1,x_1,\dots,x_k)$ and
$p=1+k$.  Then $P(0)$ is known as the \emph{correlation polytope}, a
well studied polytope in combinatorial optimization.  This case is
particularly well-behaved, well-studied, and relevant for
practitioners.  It was investigated in depth in
\cite{grasshoff2013optimal,grasshoff2014optimal,grasshoff2015stochastic,kahle2015algebraic}.

The \emph{pairwise interaction model} arises for $d=2$, where
\[f(x) = (1,x_1,\dots,x_k,x_1x_2,x_1x_3,\dots,x_{k-1}x_k)\]
and $p = 1+k+\binom{k}{2}$.  This situation is already so intricate
that neither an algebraic description of the model (the set of vectors
$(\lambda (x,\beta))_{x\in\{0,1\}^k}$ parametrized by $\beta\in\RR^p$)
nor an explicit description of the polytope $P(\beta)$ are known.

An \emph{approximate design} is a vector
$(w_x)_{x\in\{0,1\}^k} \in [0,1]^{2^k}$ of non-negative weights with
$\sum_x w_x = 1$.  To each approximate design there is a matrix
$M(w,\beta) = \sum_x w_x M(x,\beta) \in P(\beta)$.  The main problem
of classical design theory is to find designs $w$ that are optimal
with regard to some criterion.  We limit ourselves to $D$-optimality,
where the determinant ought to be maximized.  To simplify the problem,
we also only consider maximizing the determinant over $P(\beta)$, and
not finding explicit weights $w$ that realize an optimal matrix
in~$P(\beta)$.  In non-linear regression, such as the Poisson
regression considered here, this optimal solution depends on~$\beta$
(in linear regression it does not).  Our approach is to consider the
set of optimization problems for all $\beta$ and subdivide them into
regions where the optima are structurally similar.  These regions of
optimality are semi-algebraic.

In our setting, there are always matrices with positive determinant in
$P(\beta)$.  Since the vertices are rank one matrices, the optimum
cannot be attained on any face that is the convex hull of fewer than
$p$ vertices.  A design $w$ is \emph{saturated} if it achieves this
lower bound, that is, $|\supp (w)| = p$.

As the logarithm of the determinant is concave, for each given
$\beta$, the optimization problem can be treated with the tools of
convex optimization.  The design problem is to determine the changes
in the optimal solution as $\beta$ varies.
% For example, an easy observation is that when the $\beta_\emptyset$
% changes, the determinant of $M(w,\beta)$ is globally scaled.  For
% all question regarding optimal design we may therefore assume
% $\beta_\emptyset = 0$.

A special design, relevant for practitioners and studied
in~\cite{kahle2015algebraic}, is the \emph{corner design $w_{k,d}^*$}.
It is the saturated design with equal weights $w_{x}=1/p$ for all
$x \in \{0,1\}^k$ with $|x|_1 \le d$.  For example, for $k=3$ rules
and interaction order $d=2$ the regression function is
$f(x_1,x_2,x_3)=(1,x_1,x_2,x_3,x_1x_2, x_1x_3, x_2x_3)$ and there are
$p=7$ parameters.  The corner design has weight $1/7$ on the seven
binary 3-vectors different from~$(1,1,1)$.

Saturated designs are mathematically attractive due to their
combinatorial nature.  It is reflected in the following classical
theorem of Kiefer and Wolfowitz which is a main tool in the theory of
optimal designs.  See~\cite[Section~9.4]{pukelsheim1993optimal}
or~\cite{kiefer1960equivalence} for details and proofs.

\begin{theorem}\label{t:kw}
Let $X\subset \{0,1\}^k$ be of size~$p$.  There is a matrix with
optimal determinant in the face $\conv\{M(x,\beta) : x\in X\}$ if and
only if for all $x\in\{0,1\}^k$
\[
\lambda(x,\beta) (F^{-T}f(x))^T \psi^{-1}(\beta) (F^{-T}f(x)) \le 1.
\]
where $F$ is the $(p\times p)$-matrix with rows $f(x), x\in X$ and
$\psi$ is the diagonal matrix $\diag(e^{\beta_1},\dots,e^{\beta_p})$.
If this is the case, then the optimal point is
$\frac{1}{p}\sum_{x\in X} M(x,\beta)$, the geometric center of the
face.
\end{theorem}

After changing the scale by the introduction of parameters
$\mu_i = e^{\beta_i}$, Theorem~\ref{t:kw} yields a system of rational
polynomial inequalities in the~$\mu_i$.  Together with the
requirements $\mu_i > 0$, we find a semi-algebraic characterization of
regions of optimality for saturated designs.

For example, the inequalities corresponding to the corner design are
the topic of \cite{kahle2015algebraic}.  It can be seen that there
always exist parameters $\beta_1,\dots,\beta_p$ that satisfy the
inequalities in Theorem~\ref{t:kw}.  A good benchmark for our
understanding of the semi-algebraic geometry of the Rasch Poisson
counts model is to understand the other saturated designs, raised as
\cite[Question~3.7]{kahle2015algebraic}.
\begin{question}\label{q:cornerOnly}
When $\beta_i < 0$, for all $i=1,\dots,p$, is the corner design the
only saturated design $w$ that admits parameters $\beta$ such $w$ is
$D$-optimal for~$\beta$?
\end{question}
For $d=1, k=3$, Question~\ref{q:cornerOnly} has been answered by
Graßhoff et al. They have shown that, up to fractional factorial
designs at $\beta=0$, only the corner design yields a feasible
system~\cite{grasshoff2015stochastic}.  Using computer algebra, the
case $d=1, k=4$ can be attacked.  
% We describe this approach in the next section.

\section{Non-optimality of saturated designs for four predictors}

Our benchmark problem for computational treatment of inequality
systems is an extension of the content of
\cite{grasshoff2015stochastic} to the case~$d=1$ and $k=4$.  Together
with Philipp Meissner, at the time of writing a master student, we
have undertaken computational experiments.  In this situation $p=5$
and a saturated design is specified by a choice of its support
$X \subset \{0,1\}^4$ with $|X|=5$.  A number of reductions applies.
For example, if all 5 points lie in a three-dimensional cube, the
determinant can be seen to be equal to zero throughout the face, so
that optimality is precluded from the beginning.  The hyperoctahedral
symmetry acts on the designs and the inequalities.  Therefore only one
representative of each orbit has to be considered.  After these
reductions we are left with 17 systems of inequalities, one for each
orbit of supports of saturated designs.  One orbit corresponds to the
corner design for which there always exist parameters at which it is
optimal.  It is conjectured that the remaining 16 saturated designs
admit no parameters under which they are optimal.  Theorem~\ref{t:kw}
translates this conjecture into the infeasibility of 16 inequality
systems.  The most complicated looking among them is the following.
\begin{gather*}
4 \mu_1\mu_2\mu_3 \mu_4 + \mu_1 \mu_3 + \mu_1 \mu_2 + 4 \mu_2 \mu_3 +
\mu_4 - 9 \mu_2\mu_3 \mu_4 \le 0  \\
4 \mu_1 \mu_2 \mu_3 \mu_4 + \mu_2 \mu_3 +
\mu_1 \mu_2 + 4 \mu_1 \mu_3 + \mu_4 - 9 \mu_1 \mu_3 \mu_4 \le 0 \\
4 \mu_1 \mu_2 \mu_3 \mu_4 + \mu_2 \mu_3 + \mu_1 \mu_3 + 4 \mu_1 \mu_2
+ \mu_4 - 9 \mu_1 \mu_2 \mu_4 \le 0 \\
\mu_1 \mu_2 \mu_3 \mu_4 + \mu_2 \mu_3 + \mu_1 \mu_3 + \mu_1 \mu_2 + \mu_4 - 9 \mu_1 \mu_2 \mu_3 \le 0 \\
\mu_1 \mu_2 \mu_3 \mu_4 + \mu_1 \mu_3 + \mu_2 \mu_3 + 4 \mu_1 \mu_2 +
4 \mu_4 - 9 \mu_3 \mu_4 \le 0 \\
\mu_1 \mu_2 \mu_3 \mu_4 + \mu_1 \mu_2 + 4 \mu_1 \mu_3 + \mu_2 \mu_3 +
4 \mu_4 - 9 \mu_2 \mu_4 \le 0 \\
\mu_1 \mu_2 \mu_3 \mu_4 + \mu_1 \mu_2 + 4 \mu_2 \mu_3 + \mu_1 \mu_3 +
4 \mu_4 - 9 \mu_1 \mu_4 \le 0 \\
\mu_1 \mu_2 \mu_3 \mu_4 + 4 \mu_1 \mu_3 + 4 \mu_2 \mu_3 + \mu_1 \mu_2
+ \mu_4 - 9 \mu_3 \le 0 \\
\mu_1 \mu_2 \mu_3 \mu_4 + 4 \mu_1 \mu_2 + \mu_1 \mu_3 + 4 \mu_2 \mu_3
+ \mu_4 - 9 \mu_2 \le 0 \\
\mu_1 \mu_2 \mu_3 \mu_4 + 4 \mu_1 \mu_2 + \mu_2 \mu_3 + 4 \mu_1 \mu_3
+ \mu_4 - 9 \mu_1 \le 0 \\
4 \mu_1 \mu_2 \mu_3 \mu_4 + \mu_1 \mu_2 + \mu_1 \mu_3 + \mu_2 \mu_3 +
4 \mu_4 - 9 \le 0 \\
\mu_1 > 0, \quad \mu_2 > 0, \quad \mu_3 > 0, \quad \mu_4 > 0.
\end{gather*}
The interested reader is invited to try her favorite method of showing
infeasibility of this system.  We have first tried SDP methods.  In
the best situation, they would yield an Positivstellensatz
infeasibility certificate (maybe for a relaxation).  For this we used
\yalmip~\cite{YALMIP} together with the \mosek solver~\cite{mosek} to
set up moment relaxations.  While in general this method works and is
reasonably easy to set up, it is not applicable here as the
infeasibility of the system seems to depend on the strictness of the
inequalities $\mu_i > 0$.  Since spectrahedra are closed, the SDP
method only works with closed sets.  Tricks like introducing a new
variables which represents the inverses of the $\mu_i$ lead to
unbounded spectrahedra. Bounding these is equivalent to imposing an
arbitrary bound $\mu_i \ge \epsilon$.  With this the degrees of the
Positivstellensatz certificate for infeasibility grow (quickly) when
$\epsilon \to 0$.  In total, the numerical method can give some
intuition, but it is not feasible to yield proofs for the benchmark
problem.

Our second attempt was to use \qepcad~\cite{brown2003qepcad}, a
somewhat dated open source implementation of quantifier elimination.
The system is very easy to use, but unfortunately it seems to have
problems already with small polynomial inequality systems due to a
faulty memory management in the underlying library \textsc{saclib}.
There have been attempts to rectify the
situation~\cite{Richardson2005}, but their source code is unavailable
and the authors are unreachable.

Finally, we tried the closed source implementation of quantifier
elimination in \mathematica~\cite{mathematica} and were positively
surprised about its power.  Its function \textsc{Reduce} quickly
yields that \textsc{FALSE} is equivalent to the existence of
$\mu_1,\dots,\mu_4$ satisfying some of the 17 inequality systems.
However, the benchmark system above seems out of reach.  From here,
the road is open to trying various semi-automatic tricks.  For
example, \mathematica can confirm within a reasonable time frame that
there is no solution to the above inequality system when
$\mu_3 = \mu_4$ is also imposed.  A summary of our findings will
appear in the forthcoming master thesis of Philipp Meissner.
% of The computation takes less than a second on a standard laptop.
% The other 16 systems can be studied similarly.  If one accepts
% closed source computational results, a generalization of the results
% of \cite{grasshoff2015stochastic} to the case $k=4$ is thereby
% proven.

\section{Outlook}

Whoever takes an experimental stance towards mathematics will, from
time to time, be faced with polynomial systems of equations and
inequalities.  We have shown one such a situation coming from
statistics here and there are more to be found from the various
equivalence theorems in design theory~\cite{pukelsheim1993optimal}.

Deciding if such a system has a solution is a basic task.  The
technology to solve it should be developed to a degree that a
practitioner can just work with off the shelf software to study their
polynomial systems.  For systems of equations this is a reality.
There are several active open source systems that abstract Gröbner
bases computations to a degree that one can simply work with
ideals~\cite{CoCoA-5,M2,Singular}.  For systems of polynomial
inequalities, the situation is not so nice.  The method to exactly
decide feasibility of general polynomial inequality systems is
quantifier elimination~\cite[Chapter~14]{basu2006algorithms}.  The
only viable open source software for quantifier elimination is \qepcad
which appears unmaintained for about a decade.  There do exist closed
implementations that seem to work much better, for example in
\mathematica.  Whether one accepts a proof by computation in a closed
source system is a contentious matter.

\begin{problem}
Develop a fast and user-friendly open source tool to study the
feasibility of polynomial inequality systems with quantifier
elimination.
\end{problem}

We shall not fear the complexity theory.  The documentation and use
cases of \qepcad demonstrate that many interesting applications were
in the reach of quantifier elimination already a decade ago.  Gröbner
bases were deemed impractical in view of their complexity theory, yet
they are an indispensable tool now.  We hope that in the future exact
methods in semi-algebraic geometry can be developed to the same extend
as exact methods in algebraic geometry are developed.  

Finally, for experimentation one can always resort to numerical
methods.  Via the Nullstellensatz and the various Positivstellensätze
the optimization community has developed very efficient methods to
deal with polynomial systems of equations and
inequalities~\cite{de2012computation}.

% \begin{problem}
% Use general equivalence theorems to derive semi-algebraic
% characterizations of regions of optimality in other non-linear
% regression models.
% \end{problem}

% \bibliographystyle{amsplain}
% \bibliography{poisson}

\providecommand{\bysame}{\leavevmode\hbox to3em{\hrulefill}\thinspace}
\providecommand{\MR}{\relax\ifhmode\unskip\space\fi MR }
\providecommand{\MRhref}[2]{%
  \href{http://www.ams.org/mathscinet-getitem?mr=#1}{#2}
}
\providecommand{\href}[2]{#2}

\end{document}